\documentclass[11pt]{article}
\usepackage{latexsym,amsmath,amssymb, amsthm}

\newif\ifdviwin

\usepackage[latin1]{inputenc}
\usepackage[english]{babel}
\usepackage{indentfirst}
\usepackage[mathscr]{eucal}
\usepackage{amssymb,amsmath,amsfonts}
\usepackage{graphicx}
\usepackage{graphicx}
\usepackage{epstopdf}
\usepackage{wrapfig}
\usepackage{float}

\newif\ifdviwin

\dviwintrue

\def\cC{\mathcal{C}}

\let\8=\infty \let\0=\emptyset

\def\flecha{\longrightarrow}

\def\cte.{\mathop{\rm cte.}\nolimits}

\def\det{\mathop{\rm det}\nolimits}

\def\re{\mathop{\rm Re }\nolimits} \def\Im{\mathop{\rm Im }\nolimits}

\def\Q{\mathbb{Q}}
\def\r{\mathbb{R}}

\def\c{\mathbb{C}}

 \newtheorem{defi}{Definition}[section] 
 \newtheorem{teo}[defi]{Theorem}

\DeclareGraphicsExtensions{.png,.pdf,.pdf,.mps,.gif,.bmp}

\begin{document}
\mbox{}\vspace{0.4cm}\mbox{}

\begin{center}
\rule{14.2cm}{1.5pt} \vspace{0.5cm}

{\Large \bf Improper affine spheres and the Hessian one equation \footnote{Research partially supported by Ministerio de Educaci\'on y Ciencia Grant No. MTM2016-80313-P and Junta de Andaluc\'{\i}a Grant No. FQM-325.}} \\
\vspace{0.5cm}
{\large Antonio Mart\'{\i}nez and Francisco Mil\'an}\\
\vspace{0.3cm} \rule{14.2cm}{1.5pt}
\end{center}
Departamento de Geometr\'{\i}a y Topolog\'{\i}a, Universidad de Granada, E-18071 Granada, Spain. \\
e-mail: amartine@ugr.es, milan@ugr.es
\vspace{0.2cm}
\\
{\sl Keywords:} Hessian one equation, Improper affine spheres, Ribaucour transformations, Cauchy problem.
\\
{\sl 2000 Mathematics Subject Classification:} 53A15.
\vspace{0.3cm}

\begin{abstract}
Improper affine spheres have played  an important role in the development of geometric methods for the study of the Hessian one equation. Here, we  review most of the advances we have made in this direction during the last twenty years. \end{abstract}
\section{Introduction}

Differential Geometry and Partial Differential Equations (PDEs) are related by a productive tie by means of which both theories well out. 
On the one hand, many results from the theory of hypersurfaces and submanifolds use, in a fundamental way, tools from PDEs theory. Among them, we can emphasize topics like maximum
principles, regularity theorems, height estimates, asymptotic behaviour at infinity, representation
theorems, or results of existence and uniqueness with fixed boundary conditions.

\

On the other hand, many classic PDEs are linked to interesting geometric problems. Even more, apart
from being source of inspiration in the search for interesting PDEs, the geometry allows, in many cases,
to integrate these equations, to establish non trivial properties of the solutions and to give superposition
principles which determine new solutions in terms of already known solutions.

One of the biggest contributions from geometry to PDEs theory are Monge-Ampère equations. Such
equations are totally non linear PDEs which model interesting geometric aspects related to the
curvature, and its study has become a topic of great mathematical importance, see  \cite{LJSX, TW}. 

\

Among the most
outstanding Monge-Ampère equations we can quote the  Hessian one equation 
\renewcommand{\theequation}{$H_\epsilon$}
\begin{equation} 
u_{xx} u_{yy} - u_{xy}^2=\varepsilon, \qquad \varepsilon\in \{-1,1\}, \qquad (x,y)\in \Omega\subseteq \mathbb{R}^2  .\label{hess}
\end{equation}
The equation \eqref{hess} has been studied from many perspectives by several authors and the situation changes completely if we take $\varepsilon=1$ (elliptic case) or $\varepsilon=-1$ (hyperbolic case). When $\varepsilon=1$,
we know there exists   a representation for the solutions in terms of holomorphic data  and their
most relevant global properties has been studied, such as the behaviour of the ends, the flexibility or the
existence of symmetric solutions (see  \cite{FMM1, FMM2}). In \cite{CL, FMM3} the space of solutions of the exterior Dirichlet problem is endowed with a suitable structure and in  \cite{GMMi}, the space of solutions of the equation defined
on $\mathbb{R}^2$ with a finite number of points removed was described, and endowed with the structure of a differential manifold. In \cite{ACG} the Cauchy problem was solved, and this solution was
used, among other things, to classify the solutions which have an isolated singularity. Thanks to
classical subjects of the PDEs theory, such as the continuity method, it seems possible to extend some
of the above results to more general classes of Monge-Ampère equations.  However, the hyperbolic  case is much more complicated and we can not expect  classification results as in the definite case. For example, $u(x,y)= xy + g(x)$ is an entire solution for any real function $g$.

\renewcommand{\theequation}{\arabic{equation}}

\

 From a geometric point of view, the equation \eqref{hess}  arises as the equation of  improper affine spheres, that is,  surfaces with
parallel affine normals.  Although Affine Differential Geometry has a long history whose origins date back to 1841 in a work by Transon on the normal affine of a curve, it has been during the last thirty years when this theory has experienced a remarkable development. In its research geometric, analytic and complex techniques
mix and in order to understand the geometry of this important family of surfaces
deeply, the global rigidity has been weakened in different ways. On the one hand, in \cite{FMM1,FMM2}, using methods of complex analysis, their behaviour at the infinity has been described  noting that there exists a closely relationship between their ends and the ones of a minimal surface with finite total curvature. On the other hand and 
concerning with \eqref{hess}, a geometric theory of smooth maps with singularities (improper affine maps) has been developed  opening an interesting range of global examples. In most of the cases the singular set determines the surface and, generically, the singularities are cuspidal edges and swallowtails, see \cite{CST,IM,M}.

\

In very recent works, \cite{MM2, MMT, M1, M2} we have solved the problem of finding all indefinite improper affine spheres passing through a given regular curve of $\r^3$ with a prescribed affine co-normal vector field along this curve; the problem is well-posed when the initial data are non-characteristic and the uniqueness of the solution  can fail at characteristic directions.  We also  have learnt how to obtain easily improper affine maps with a prescribed singular set, how to construct  global examples  with the desired singularities and how to use the  classical theory of Ribaucour transformations   to obtain new solutions of the elliptic hessian one equation. 

\

Here, we want to make a short survey with some of the  above mentioned advances in the geometric study of  \eqref{hess}.  
\section{A  complex (split-complex) representation formula}
If $u : \Omega \subset \r^2 \longrightarrow \r^3$ is a solution of \eqref{hess}, 
then its graph $$\psi = \{(x, y, u(x,y)) : (x, y) \in \Omega \}$$ is an  improper affine sphere in $\r^3$
with constant affine normal $\xi=(0,0,1)$, affine metric $h$, 
\begin{equation}
h:= u_{xx}\ dx^2+  u_{yy}\ dy^2 + 2  u_{xy} \ dx dy, 
 \label{ametric}
\end{equation}
and affine conormal $N$,
\begin{equation}
N:= (-u_x,-u_y,1).  
 \label{conormal}
\end{equation}
From \eqref{ametric} and \eqref{conormal} it is easy to check that the following relations hold,
\begin{align}
&h = -<dN,d\psi>, \quad <N,\xi> = 1,  \quad <N,d\psi> = 0, \label{feq}\\
&\sqrt{|\det(h)|} = \det[\psi_x,\psi_y,\xi] = - \det[N_x,N_y,N], \label{metric}
\end{align}
 see \cite{LSZ,NS} for more details. Conversely, up to unimodular transformations,  any improper affine sphere in $\r^3$ is, locally, the graph over a domain in the $x,y$-plane of a solution to \eqref{hess}.

\

In the elliptic (hyperbolic) case the affine metric $h$ induces a Riemann (Lorentz) surface structure on $\Omega$ known as  the {\sl underlying conformal structure of $u(x,y)$. }

It follows from \eqref{hess} that,
\begin{equation}
(d u_x)^2 + \varepsilon \; dy^2 = u_{xx} h, \qquad (d u_y)^2 + \varepsilon \; dx^2 = u_{yy} h,\label{conf}
\end{equation} 
and the two first coordinates of $\psi$ and $N$ provide conformal parameters for $h$. Actually, consider  $\c_\varepsilon $ the complex (split-complex) numbers, that is 
\begin{equation}
\c_\varepsilon  = \{ z = s + j\ t \ : \ s,t\in \r, \  j^2 = -\varepsilon, \ j 1 = 1 j \}, \label{cn}
\end{equation} 
then it is not difficult to  prove, see \cite{C2,FMM2,M1}, that  $\Phi:\Omega\longrightarrow \c^3_\varepsilon$, 
\begin{equation}
\Phi:= N+ j\ \xi\times u, \label{Phi}
\end{equation} is a planar holomorphic (split-holomorphic) curve. In fact, $\Phi=(-B,A,1)$ where
\begin{equation}
A:=  -u_y + j\ x, \qquad  B:= u_x + j\ y, \label{ab}
\end{equation}
are holomorphic (split-holomorphic)  functions on $\Omega$. Moreover, from \eqref{hess} and \eqref{ametric},
\begin{equation}
|d\Phi|^2 = |dA|^2 + |dB|^2 = (u_{xx} + u_{yy}) h, 
\end{equation}
and $|d\Phi|^2 $ and  $h$ are in the same conformal class always that  $u_{xx} + u_{yy}$ has a sign.

From \eqref{ametric} and \eqref{ab}, the metric $h$ is given by 
\begin{equation}
h:= \Im (dA \overline{dB}) =  | dG|^2 - |dF|^2 \label{cmetric}
\end{equation}
where $ 2 F = -B - \varepsilon j \ A $ and $ 2 G = B- \varepsilon j\ A$, and 
 the immersion  $\psi$ may be recovered as
\begin{equation} 
\psi:= -\frac{1}{2}\Im \int  (\Phi +\overline{\Phi})\times d\Phi = ( G + \overline{F},  \frac{|G|^2}{2} -  \frac{|F|^2}{2} + 2 \re \int G  d \overline{F} ), \label{ias}
\end{equation}
with  the two first coordinates of $\psi$ identified as  numbers of $\c_\varepsilon$ in the standard way. 

Conversely, if we consider the class of improper affine maps, {\rm (IAM)}, as the improper affine spheres with admissible singularities, where the affine metric is degenerated, but the affine conormal is well defined, then we have

\begin{teo}{{\rm {\bf \cite{FMM1, M}}}} Let $\Sigma$ be a Riemann surface and 
 $\psi:\Sigma \longrightarrow  \c_\varepsilon \times \r$  be an improper affine map, then $\psi$ can be represented as
\begin{equation} \psi = \Big(G + \overline{F}, \frac{1}{2} |G|^2 - \frac{1}{2} |F|^2 + Re(GF) - 2 Re \int F dG\Big)
\label{cr}
\end{equation}
with $F, G:  \Sigma\longrightarrow  \c_\varepsilon $ holomorphic (split-holomorphic)  functions. In this case, the affine conormal and the affine metric of $\psi$ are given,  respectively,  by $N = (\overline{F} - G, 1)$ and $$h = |dG|^2 - |dF|^2.$$ 
\end{teo}
The pair $(G,F)$ is called the {\sl Weierstrass data of $\psi$}.

\section{The elliptic case}
The above complex representation has been an  important tool in the description of global examples and in the methods used to understand the geometry of this important family of surfaces deeply.   
Classical examples can be described as follows:
\begin{enumerate}
\item \emph{The elliptic paraboloid:} It can be obtained by taking,
 $\Sigma = \c$ and  Weierstrass data $(z,k z)$, where $k$
  is constant. (see Figure \ref{fi1}). 
\item \emph{Rotational improper affine maps:} They are obtained by taking
$\Sigma = \c\backslash \{0\}$
and  Weierstrass data
 $(z,\pm r^2/z)$, $r\in \r\backslash \{0\}$, (see Figure
\ref{fi1}).
\begin{figure}[H]
\begin{center}
\DeclareGraphicsExtensions{eps}
\includegraphics[width=.3\textwidth]{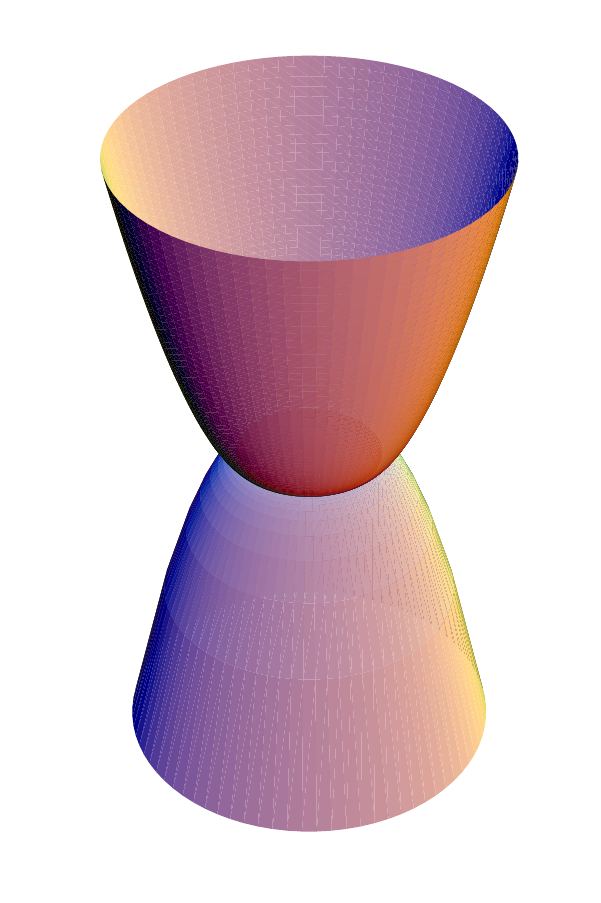}
\includegraphics[width=.3\textwidth]{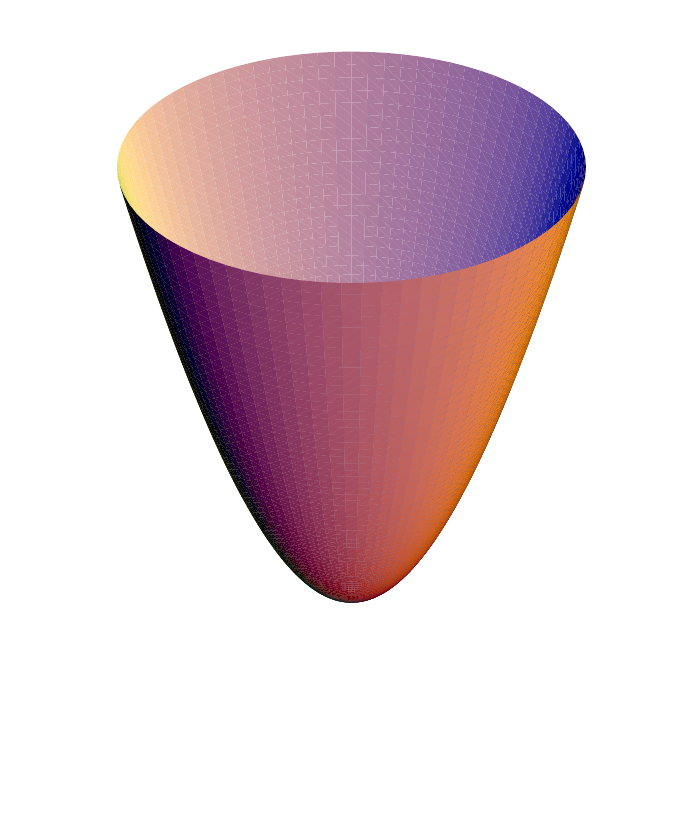}
\includegraphics[width=.25\textwidth]{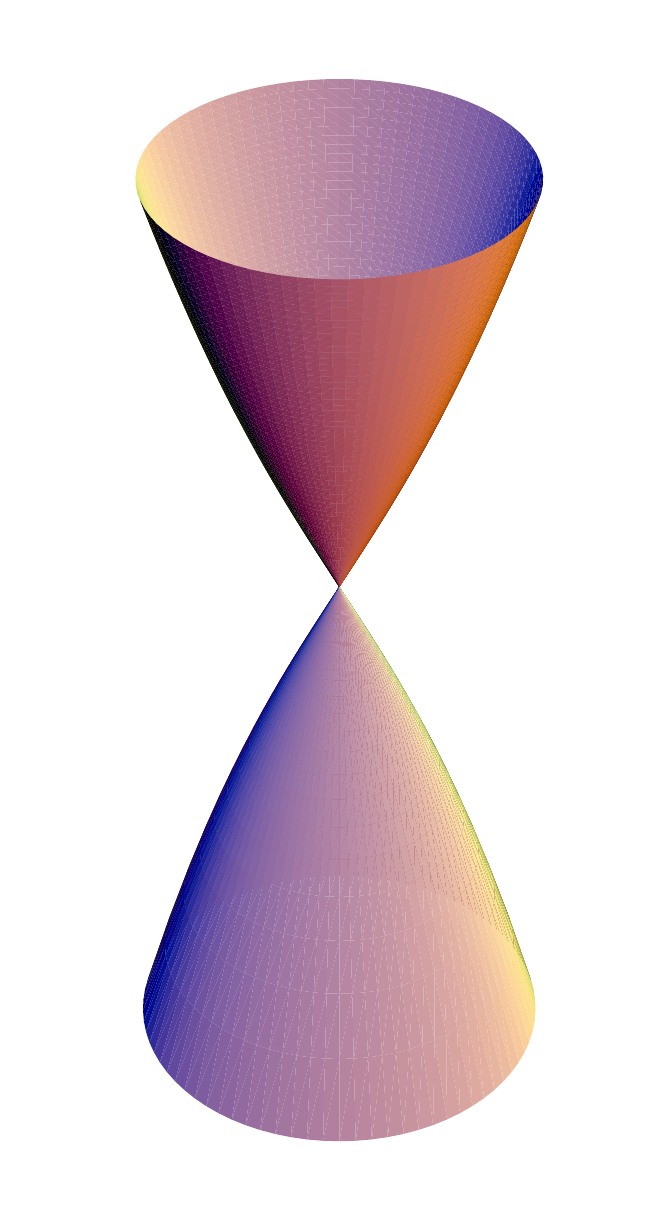}\\
\caption{Elliptic paraboloid and rotational improper affine maps}\label{fi1}
\end{center}
\end{figure}
\item \emph{Non-rotational improper affine maps with two-ends.} Others complete examples with 
two embedded ends can be obtained by taking $\Sigma =
\c\backslash \{0\}$ and  Weierstrass data: $(z, a z+ b/z + c)$,  where $b\in \r$,
$a, c\in\c$, $|a|\neq 1$ (see Figure \ref{fi2}).
\begin{figure}[H]
\begin{center}
\DeclareGraphicsExtensions{eps}
\includegraphics[width=.4\textwidth]{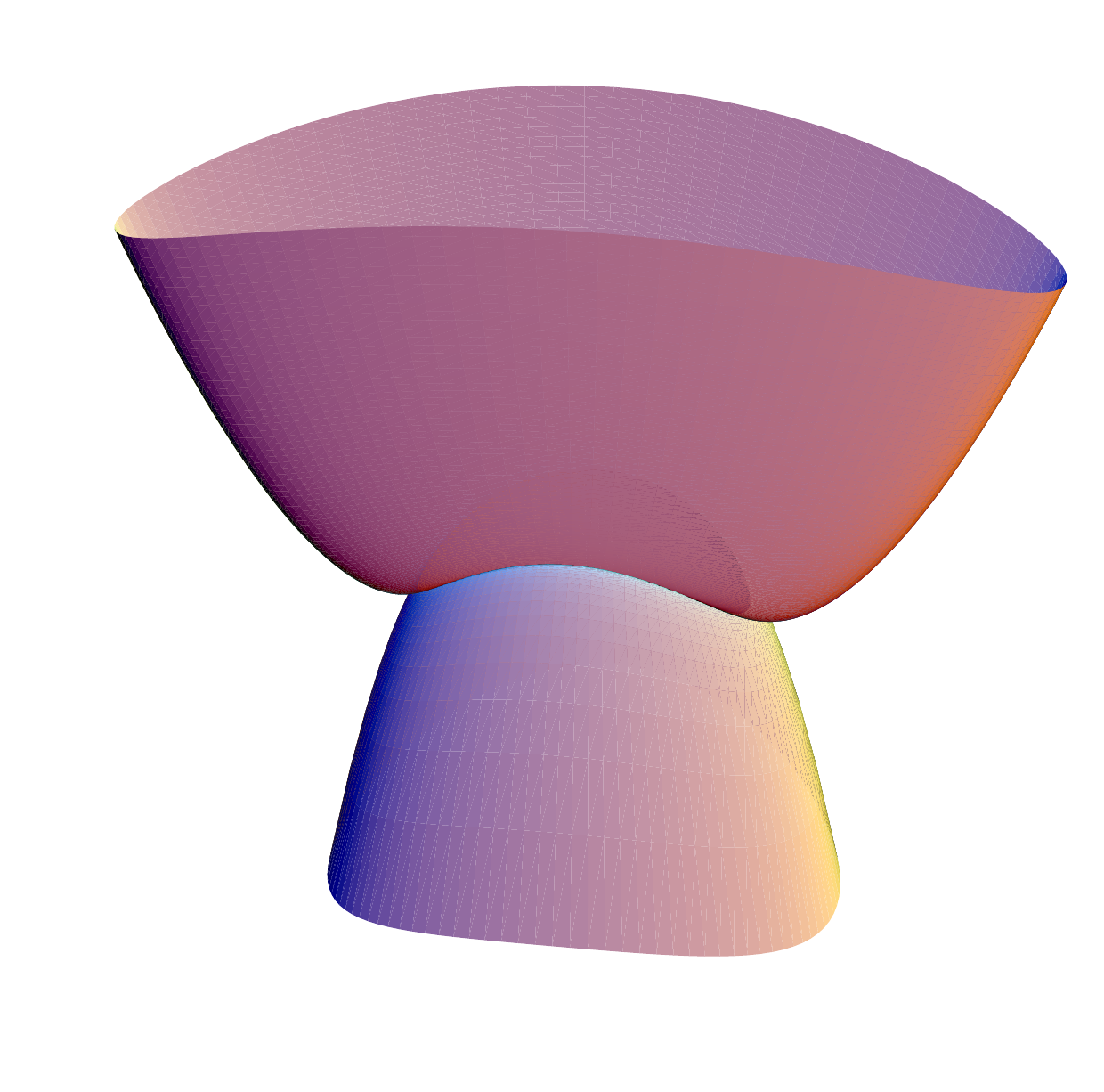}\\
\caption{Non-rotational complete improper affine map with two embedded ends}\label{fi2}
\end{center}
\end{figure}
\item \emph{Some multivalued improper affine map:} By considering $\Sigma = \c\backslash \{0\}$ and
Weierstrass data  $(z,\pm \sqrt{-1} r^2/z)$, $r\in \r\backslash \{0\}$ we obtain  multivalued  complete improper affine maps with a vertical period and two embedded ends. (see Figure \ref{fi3})
\begin{figure}[H]
  \begin{center}
\DeclareGraphicsExtensions{eps}
\includegraphics[width=.25\textwidth]{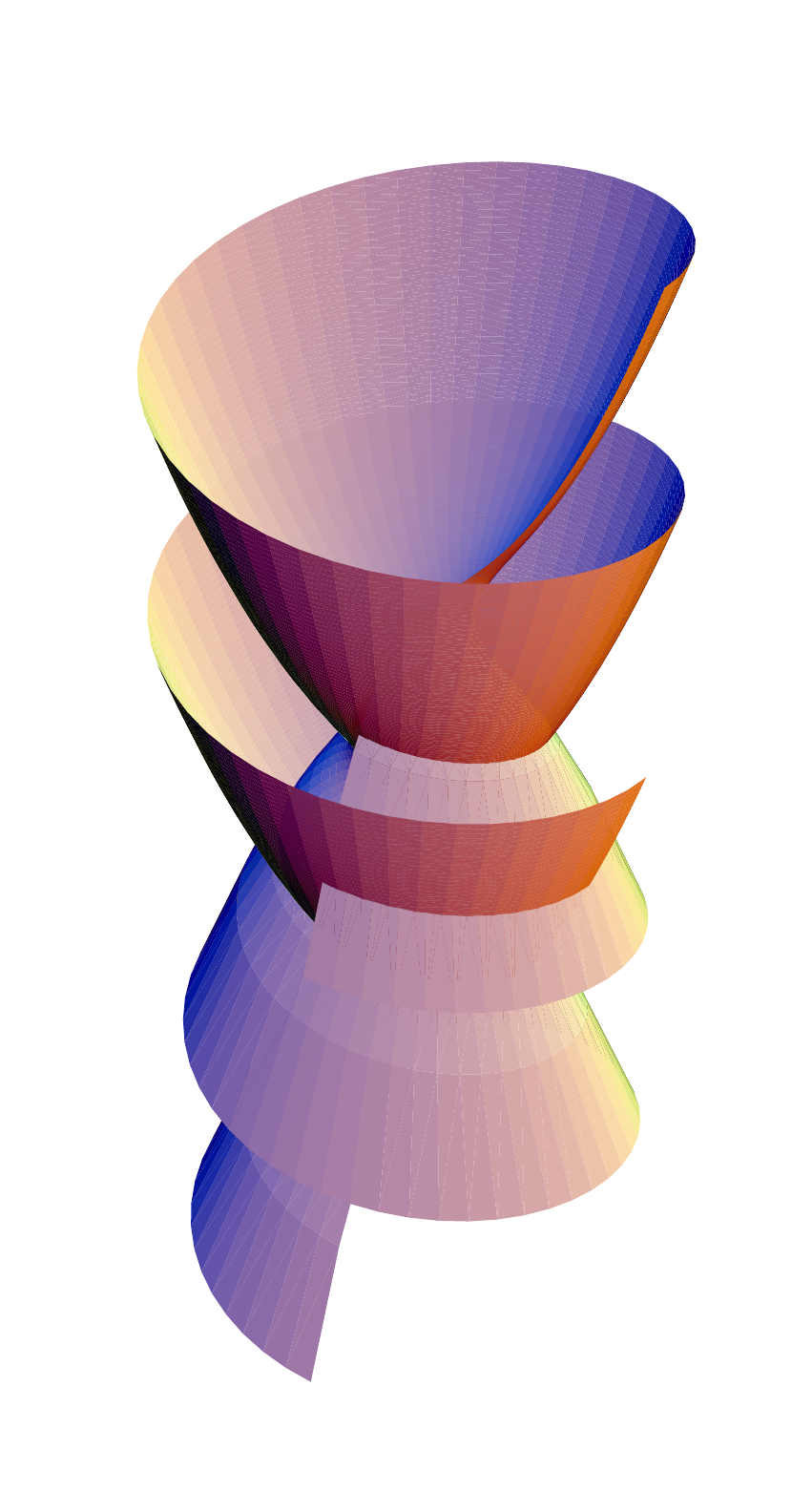}\\
\caption{Multivalued  complete improper affine map with two embedded ends }\label{fi3}
\end{center}
 \end{figure}
 \item \emph{Improper affine map with  one non-embedded end:} By taking $\Sigma = \c$ and the
Weierstrass data $(z, z + z^2 )$ we obtain a complete improper affine
 map with only one end which is non-embedded and non
 affinely equivalent to the elliptic paraboloid (see Figure \ref{fi5})
  \begin{figure}
\begin{center}
\includegraphics[width=.4\textwidth]{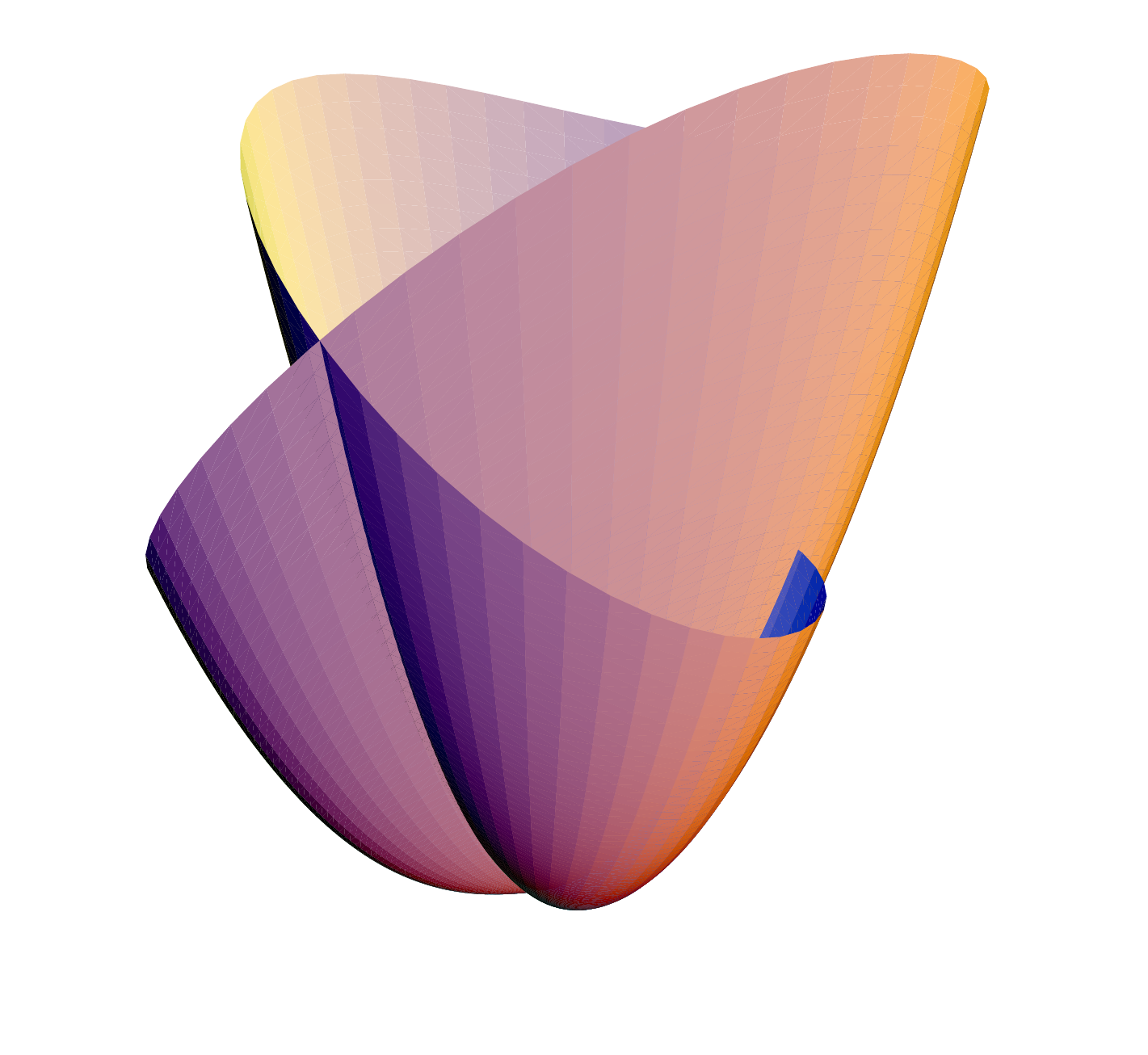}\\
\caption{Complete improper affine map with only one non-embedded end }\label{fi5}
\end{center}
\end{figure}
\end{enumerate}
\subsection{Global behaviour}
The most elemental elliptic model  is the elliptic paraboloid, it  is the
only elliptic improper affine sphere which is complete from the affine point of view, \cite{C1,CY}  and also from the
Euclidean one \cite{Jo}.

\ 

Using methods of complex analysis, \cite{M}, one can  weakened the above global rigidity and prove that any  complete IAM is conformally equivalent to a compact Riemann surface $\overline{\Sigma}$ with a finite number of points (ends) removed. In this  case the Weierstrass data extend meromorphically to the ends and 
a complete  end $p$ is embedded if and only if the Weirestrass data have at most a single pole at $p$. 
In this context, the classical result of J\"orgens is extended to prove that the elliptic paraboloid is the unique embedded complete IAM with only one end and any complete IAM with exactly two embedded ends is affine equivalent to either a rotational IAM or to one of the examples described above in 3.

These global results and \eqref{cr}, let us to prove:
\begin{itemize}
\item  Any solution of \eqref{hess} has the following  behaviour  at infinity:
$$u(x, y) \approx {\cal E}(x, y) + a \log|z|^2,$$
where ${\cal E}(x, y)$ is a quadratic polynomial, \cite{FMM2}.
\item The moduli space of the exterior Dirichlet problem associated to \eqref{hess} is either empty or a 5-dimensional differentiable manifold, \cite{FMM3}.
\end{itemize}

The conformal representation \eqref{cr} and the existence of two classical holomorphic functions with the property of mapping bijectively a once punctured bounded domain in $\c$ with $n$ boundary components onto $\c$ with $n$ vertical (resp. horizontal) slits is used in \cite{GMMi} to construct solutions 
of the elliptic Hessian one equation in the punctured plane $\r^2\backslash \{p_1,\cdots, p_n\}$ so that $p_j$'s are non removable isolated singularities of such solutions. These solutions are, up to equiaffine transformations, the only elliptic solutions to \eqref{hess} that are globally defined over a finitely punctured plane and so, there is an explicit one-to-one correspondence between the conformal equivalence classes of $\c$  with n disks removed and the space of solutions to the elliptic Hessian one equation  in a finitely punctured plane with exactly $n$ non-removable isolated singularities.

\subsection{Ribaucour's type  transformations}
It is well known  that geometric transformations are also useful  to construct new surfaces  from a given one. Recently,  the  classical study about Ribaucour transformations developed by Bianchi, \cite{Bi}, has been applied  to provided  global description of  new  families of  complete minimal surfaces obtained from minimal  surfaces which are invariant under a one-parametric group of transformations, see \cite{ CFT}. 
For the particular case of IAM we can introduce a Ribaucour's type transformation as follows:

\

Let $\psi: \Sigma \longrightarrow \c\times \r$ be an IAM, we say that an IAM $\widetilde{\psi}$ is {\sl R-associated} to $\psi$ if there is a differentiable function 
$g : \Sigma \longrightarrow \r$ such that
\begin{enumerate}
\item $(\psi + g N)  \times \xi = (\widetilde{\psi} + g \widetilde{N})  \times \xi$.
\item $dG dF = d\widetilde{G} d\widetilde{F}$, 
\end{enumerate}
where $(G,F)$ and $(\widetilde{G},\widetilde{F})$ are Weierstrass data for $\psi$ and $\widetilde{\psi}$, respectively.

\

In  \cite{MMT} we prove that $\psi$ and $\widetilde{\psi}$ are  R-associated if and only if: 
$$(\widetilde{F}, \widetilde{G}) = \Big(F + \frac{1}{cR}, G + R\Big),$$
where  $c \in \r - \{0\}$ and $R$ is a holomorphic solution of the following Riccati equation: 
\begin{equation}dR + dG = c R^2 dF \quad \Big( {\rm equivalently,} \ \  d\Big(\frac{1}{cR}\Big) + dF = \frac{1}{c R^2} dG \Big),\label{re}
\end{equation}

\begin{figure}[H]
\begin{center}
\includegraphics[width=0.16\linewidth]{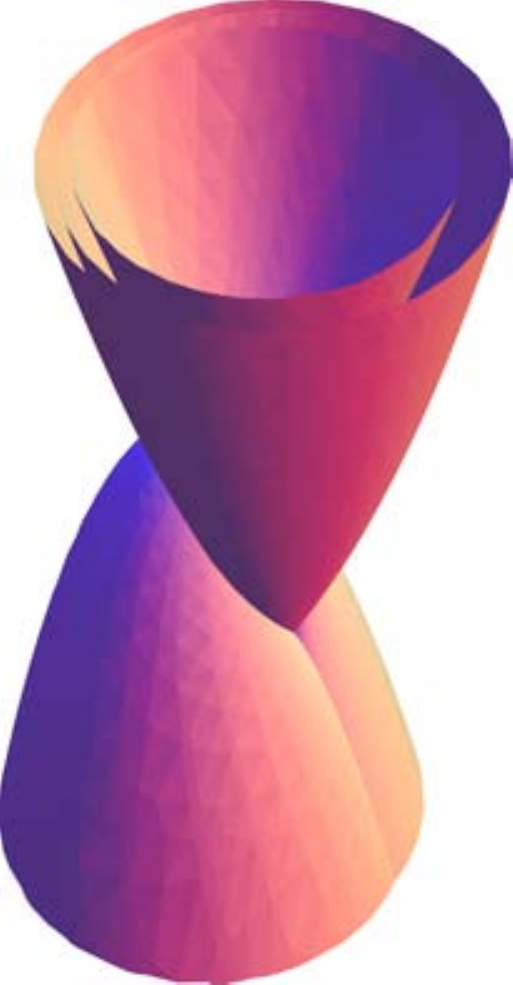}  \hspace{1cm}
\includegraphics[width=0.26\linewidth]{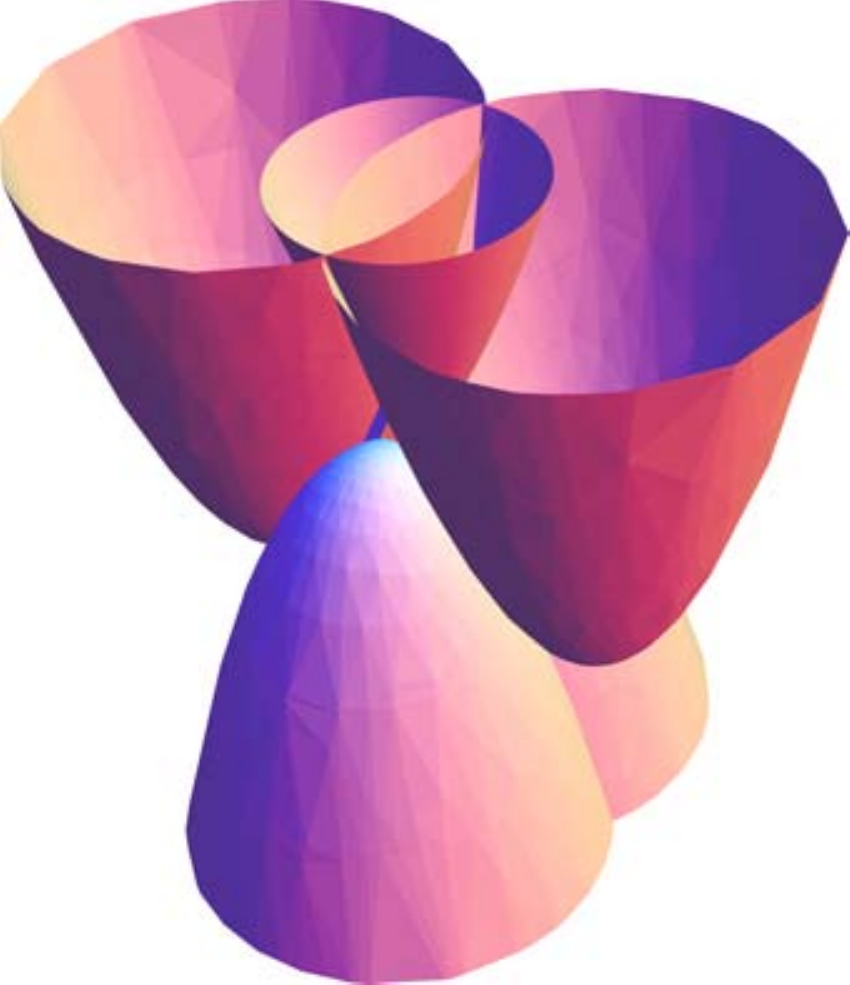}  \hspace{1cm}
\includegraphics[width=0.17\linewidth]{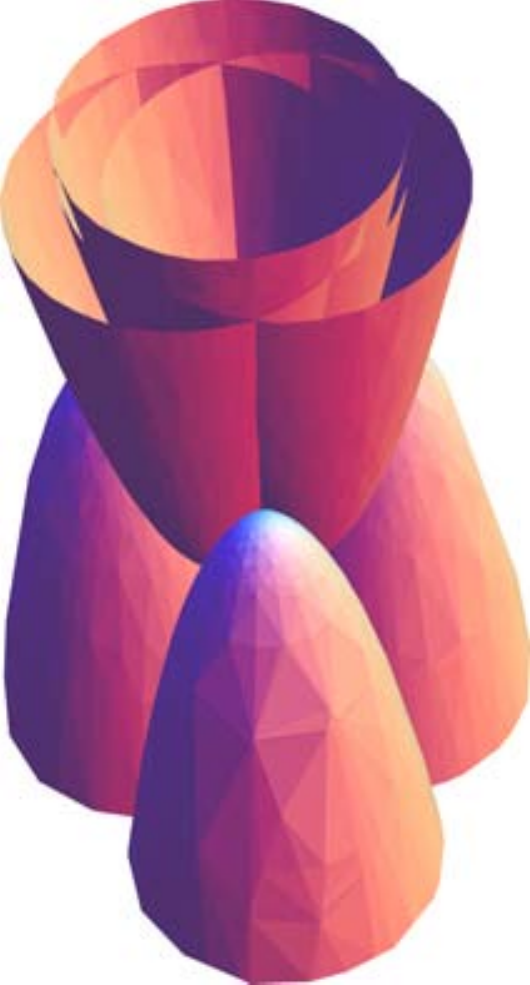}
\end{center}
\caption{ $a=1$, $n=1,\ m=3 $;  \hspace{.3cm}$a=1, \ n=2, \ m=1$;   \hspace{.3cm}   $a=1 \ n=3, \ m=2$}\label{rib2}
\end{figure}

Using \eqref{cr} and standar properties of the solutions of \eqref{re} we have  that if $\widetilde{\psi}$ is {\sl R-associated} to $\psi$, then
\begin{itemize}{\rm 
\item $\widetilde{\psi}$ has a new end at $p_0$ if and only if $p_0$ is either a zero or a pole of $R$.
\item The singular set of $\widetilde{\psi}$ is the nodal set of the harmonic function $$\log|dG| - \log(c^2|R|^4|dF|).$$ 
\item If $\psi$ is helicoidal (which means that $F G=-a^2 $ for some constant $a\in \c$) we have that the general solution of \eqref{re} can be written as  
$$R = \frac{\exp(z)}{2 a c} \frac{1 + b + (1 - b) k \exp(b z)}{1 + k \exp(b z)}$$
with  $k \in \c$,  $c \in \r - \{0\}$ and $b = \sqrt{1 + 4 a^2 c} \ne 0$.

In particular, if $$b = \frac{n}{m} \in \Q - \{0, 1\}$$ is irreducible, then $\widetilde{\psi}$ is $2 m \pi$-periodic in one variable, its singular set is contained in a compact set and has $2n$ complete embedded ends of revolution type, see Figure \ref{rib2}}.
\end{itemize}
\section{The hyperbolic case}
The central problems of hyperbolic IAM are the same as for the elliptic case. Basically this is the problem of understanding, in a global an local way, the connection between intrinsic and extrinsic affine geometry. However in contrast to the elliptic case, the results here are in many respects far from complete. For example,  we can not expect a classification result as J\"orgens' Theorem. In fact, $$ u(x,y) = x y + g(x),  \quad  (x,y)\in \r^2$$
is an entire solution of the hyperbolic Hessian one equation for any real function $g$.

\

\begin{figure}[H]
\begin{center}
\includegraphics[width=0.3\linewidth]{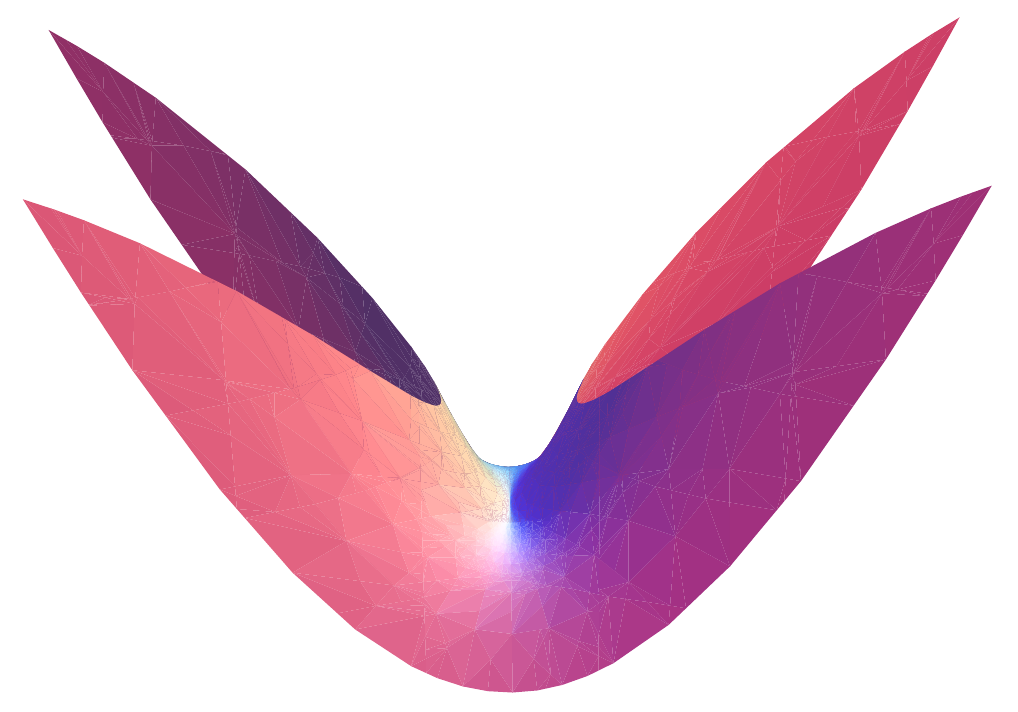}\hspace{1cm}
\end{center}
\caption{Entire solution on the puncture plane  with $H(z) = z$}
\label{fig6}      
\end{figure}

In the same way and concerning with continuous solutions of \eqref{hess}  defined on a punctured plane we have, see \cite{M2}:
\begin{itemize}
\item Let $\mathcal{D}$ be a planar disk and $A : \mathcal{D} \flecha \c_{-1}$ be a split-holomorphic function satisfying $A_z= H^2$ for some  split-holomorphic function $H : \mathcal{D} \flecha \c_{-1}$. If  $z_0 \in \mathcal{D}$ is an isolated zero of $H$, then the  indefinite improper affine map $\psi : \mathcal{D} \flecha \r^3$ given, as in \eqref{ias}, by the split-holomorphic curve $\Phi(z)=(jz, A(z),1)$,   is well defined on $\mathcal{D}^* = \mathcal{D} - \{z_0\}$ and it has a non removable  isolated singularity at $z_0$. See Fig. \ref{fig6}.
\end{itemize}
The  above result let us to construct indefinite improper affine map $\psi : \c_{-1} \flecha \r^3$ with a finite number of prescribed isolated singularities at the points $\{z_1,\cdots, z_n\}$. For this is enough to consider a split-holomorphic function $H : \c_{-1} \flecha \c_{-1}$ with zeros at the points $\{z_1,\cdots, z_n\}$.
\section{The Cauchy problem}
In many situations, the classical theory of surfaces and the methods of complex analysis have allowed  to study the classical Cauchy problem of the  underlying  associated PDE. From a geometric point of view the problem can be posed as follows:
\begin{quote}{\sc Geometric Cauchy problem:} Find the surfaces on a specific class of surfaces that pass through  a given curve  and whose tangent plane  distribution along  the curve is   also  a given distribution of oriented planes.
\end{quote}
The problem is {\sl well-posed} if there is a unique solution which depends continuously of the data. Most of the advances in this direction have been obtained in the elliptic case and for surfaces which admit a complex representation. In this case the surfaces are always real analytic and it is necessary to prescribed analytic initial data (one has  well-posedness  problems under the analytic regularity  assumption but the problem is   $\cC^\infty$ ill-posed).

In contrast to the elliptic case, the problem when the underlying PDE associated is hyperbolic, is far from complete. If one considers lesser regular data $\cC^\infty$ or even worst, many problems arise. Actually, these classes of surfaces  can not be analytic and  gluing procedures  may create unexpected situations.

In the particular case of Improper Affine Spheres (IAS), that is  improper affine maps without singularities, we want to study  the following problem:

 \begin{itemize}
 \item[({\bf P})] Find all IAS containing a curve $\alpha(s)$  with a prescribed affine conormal $U(s)$ along it.
\end{itemize}
\begin{enumerate}
\item In the elliptic case
$$0 < h(\alpha'(s), \alpha'(s)) = - \langle \alpha'(s), U'(s) \rangle,$$ 
and we will consider  $\{\alpha, U\}$ analytic data, \cite{ACG}.
\item In the hyperbolic  case, $\langle \alpha', U' \rangle$ vanishes when $\alpha'$ is an  asymptotic (also known as  characteristic) direction. In this case we only consider regular data. 
\end{enumerate}
\subsection{Non-characteristic Cauchy problem}

If $\psi:\Sigma \longrightarrow \c \times \r$ is an IAS with affine normal $\xi = (0, 0, 1)$ and $\beta : I \longrightarrow \Sigma$ is a curve, then $\alpha = \psi \circ \beta$, $U = N \circ \beta$ and $\lambda =  - \langle \alpha', U' \rangle$ satisfy
\begin{equation}
\left.
\begin{array}{l}
1 = \langle \xi, U \rangle,\\
0 = \langle \alpha', U \rangle, \\
\lambda = \langle \alpha'', U \rangle.
\end{array} \right\} \label{nonch}
\end{equation}
A pair of (analytic ) curves $\alpha, U : I \longrightarrow \r^3$ is a  {\sl non-characteritic admissible pair} if it verifies the above conditions with $\lambda : I \longrightarrow \r^+$.

\

In \cite{M1} we have proved that any non-characteristic admissible pair determines a unique IAS $\psi$ containing $\alpha(I)$ with affine conormal $U$ along $\alpha$. As consequence,  
\begin{itemize}
\item There exits a unique solution to the Cauchy problem associated to $\eqref{hess}$.

\end{itemize}

\begin{figure} [H]
\begin{center}
\includegraphics[width=.28\textwidth]{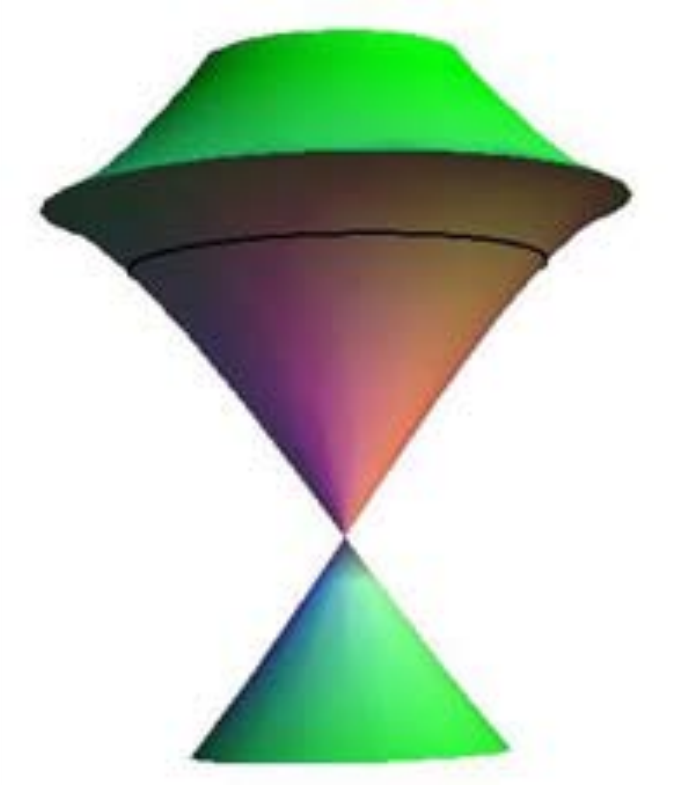} \qquad
\includegraphics[width=.4\textwidth]{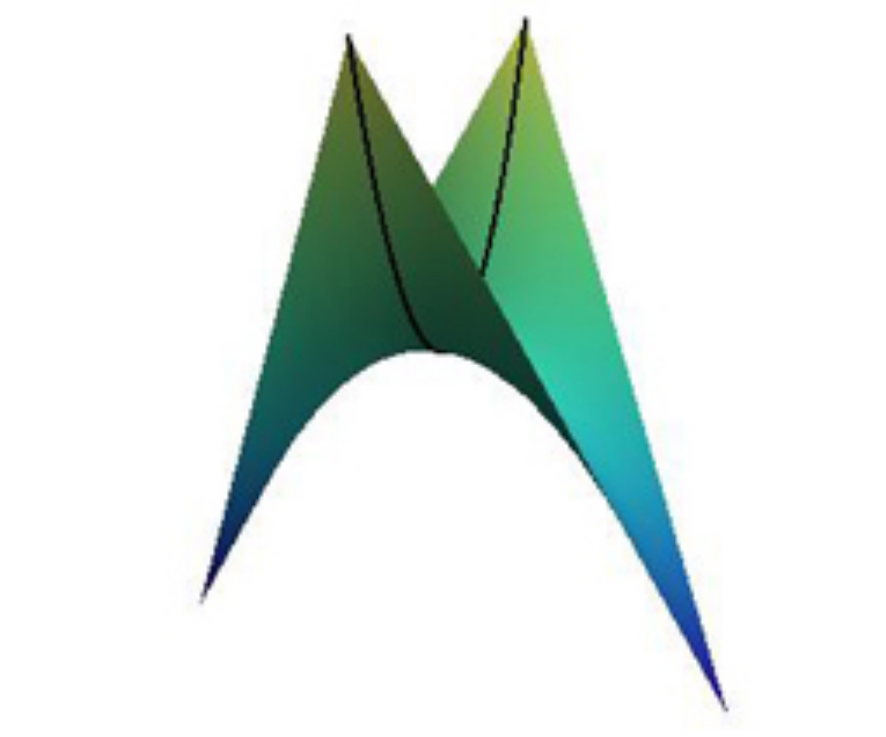}
\end{center}
\caption{IAS admitting a geodesic planar curve}\label{fi8}
\end{figure}

 If $[\alpha', \alpha'', \xi] \neq 0$, then,  from \eqref{nonch}, $U$ is determined by $\alpha$ and $\lambda$.  In particular, any  revolution IAS can be recovered with one their circles $\alpha$ and the affine 
metric along it. Moreover, $\alpha$ is geodesic when $\lambda = r^2$ and $\varepsilon = - 1$. In general, $\alpha$ is  geodesic of some IAS if and only if $[\alpha', \alpha'', \xi] = - \varepsilon [U', U'', \xi],$ with  $\lambda = m \in R^+$ (see Figure \ref{fi8}).

Any symmetry of a non-characteristic admissible pair induces a symmetry of the IAS generated by it.

\subsection{Characteristic Cauchy problem}
\begin{figure}[H] 
\begin{center}
\includegraphics[width=.4\textwidth]{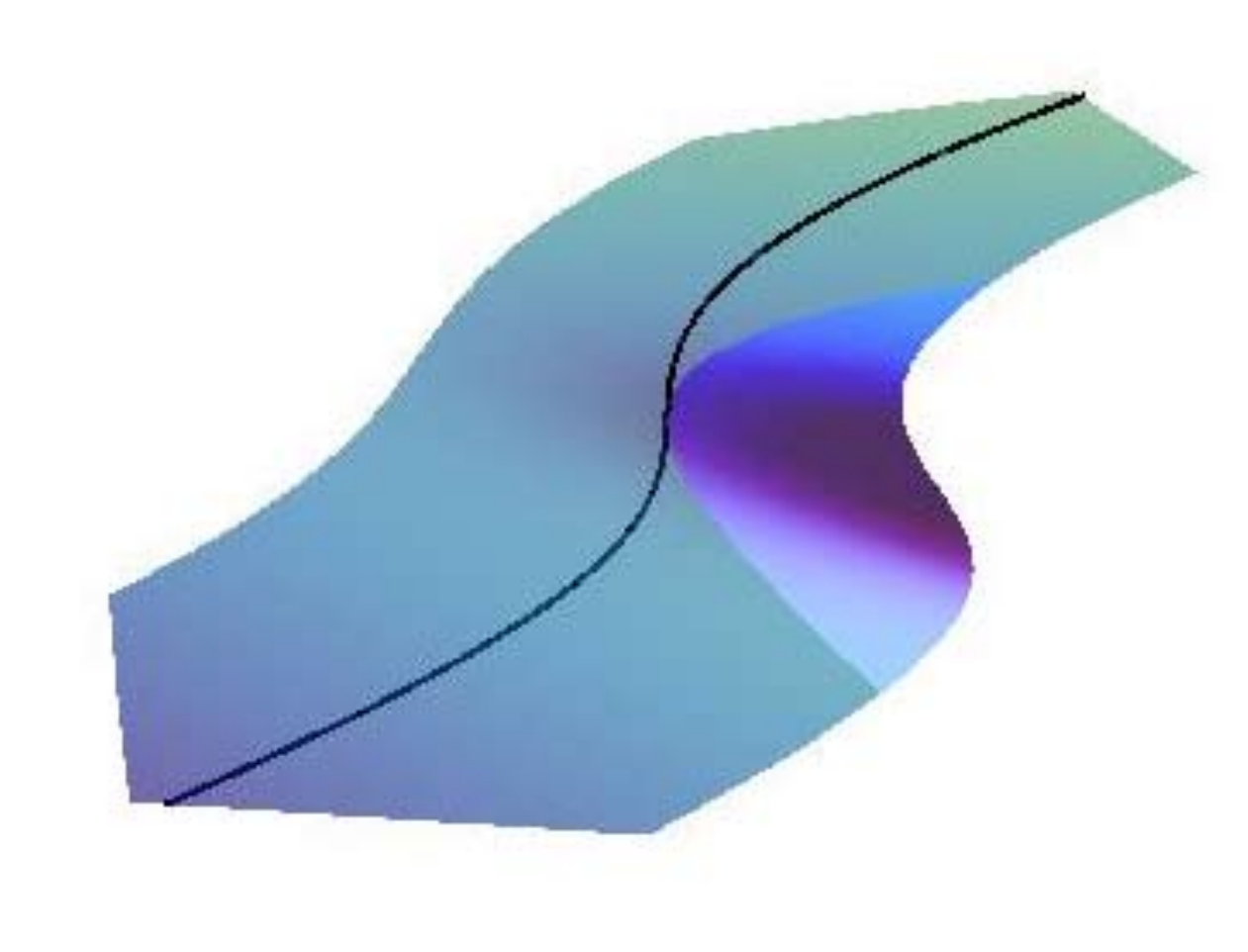} \qquad 
\includegraphics[width=.4\textwidth]{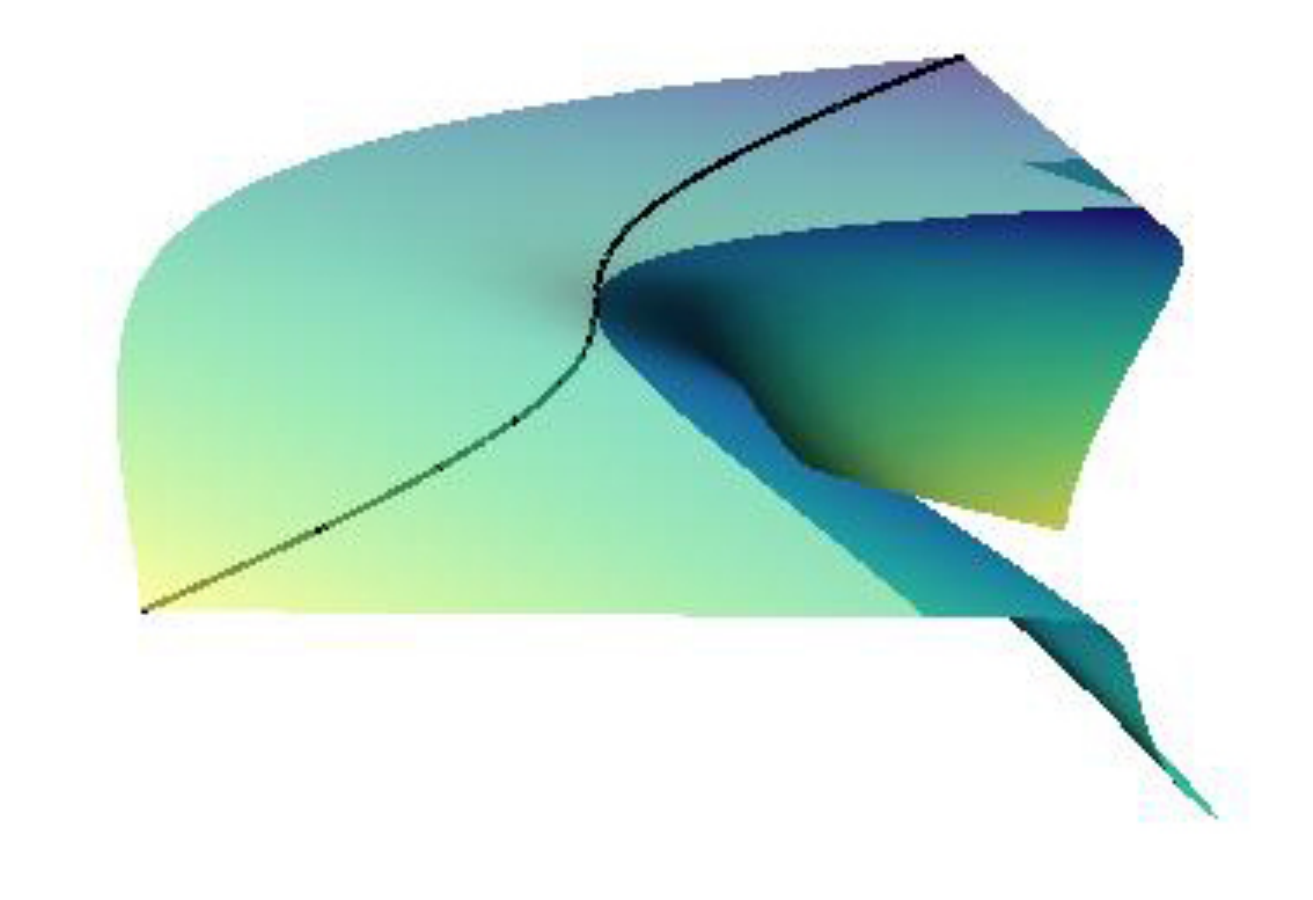}
\end{center}
\caption{Asymptotic (characteristic) data}\label{fi9}
\end{figure}

In the hyperbolic case we can consider  a characteristic admissible pair $\{\alpha, U\}$, that is, $\langle \alpha', U' \rangle$ vanishes when $\alpha'$ is an asymptotic direction. In this situation and by using  asymptotic parameters $(u, v)$, we can take,
$$N = (a + b, 1) \qquad and \qquad \xi \times \psi = (b - a, 0)$$
for  two harmonic planar curves $a(u)$ and $b(v)$ 
It is clear that an asymptotic curve $\psi(u, v_o)$ determines $a(u)$ and $N(u, v_o)$, but not $b(v)$. 
Thus, an admissible pair $\{\alpha, U\}$ generates many (hyperbolic) IAS, when $\langle \alpha', U' \rangle$ vanishes identically, in fact, the uniqueness fails when $\alpha(s) = \psi(u(s), v(s))$ is tangent to an asymptotic curve (see Figure \ref{fi9}).
\subsection{Prescribed singular curves}
\begin{figure} [H]
\begin{center}
\includegraphics[width=.25\textwidth]{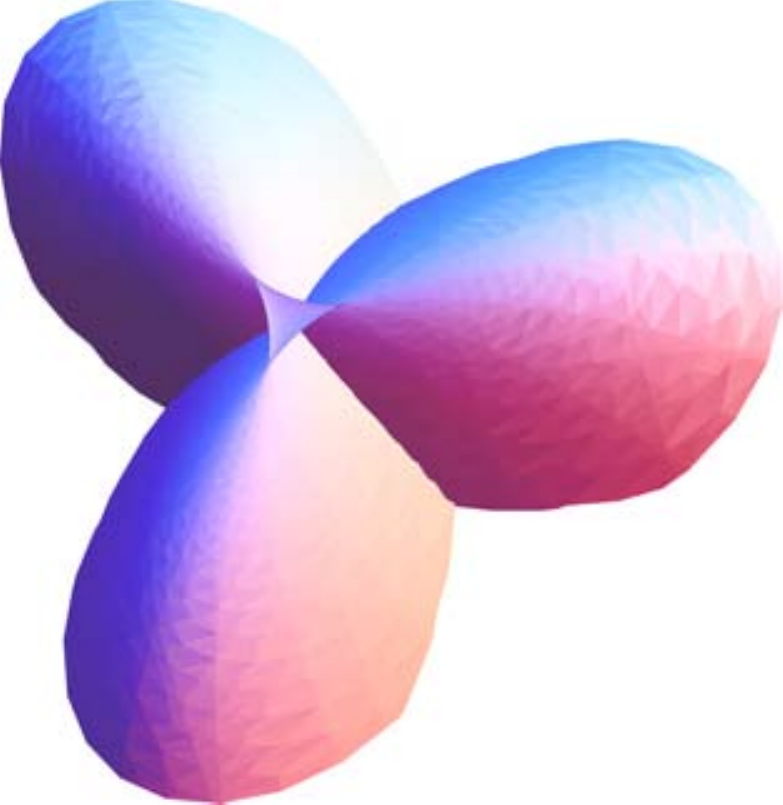} 
\includegraphics[width=.35\textwidth]{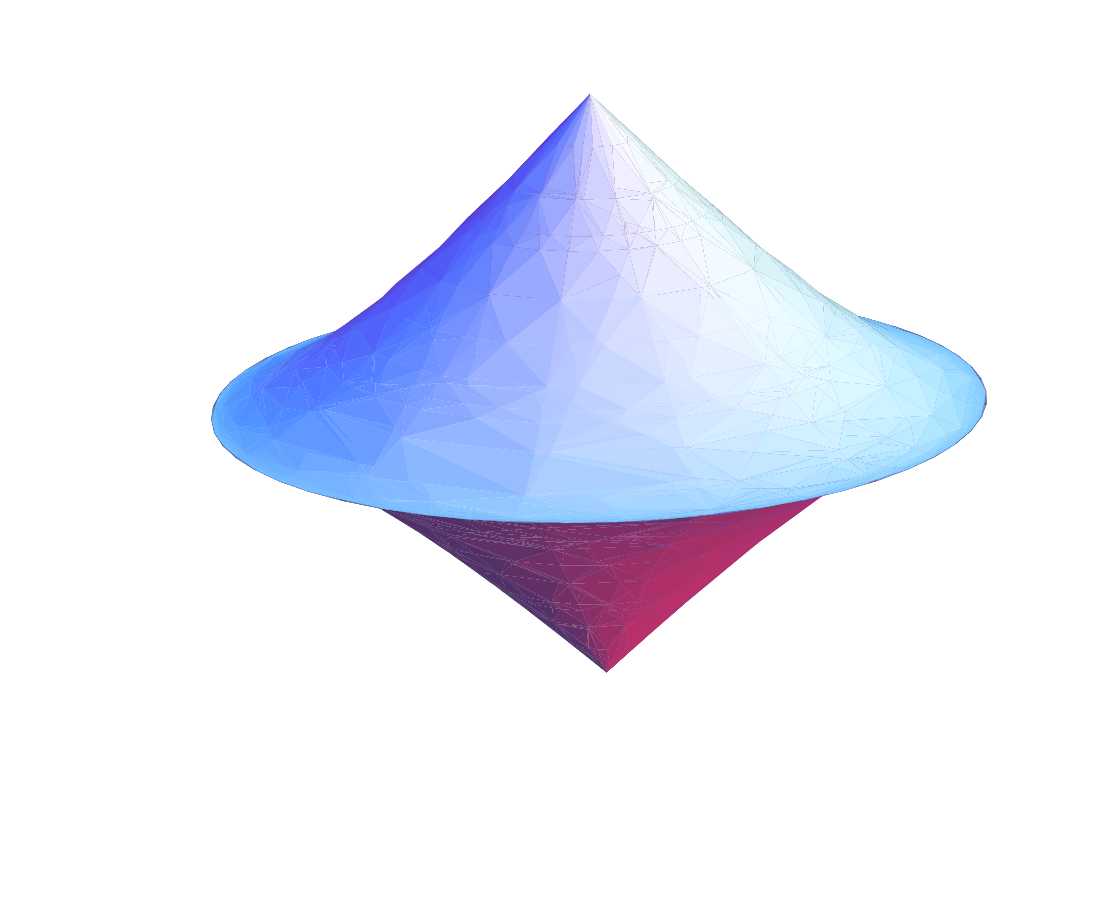}
\includegraphics[width=.28\textwidth]{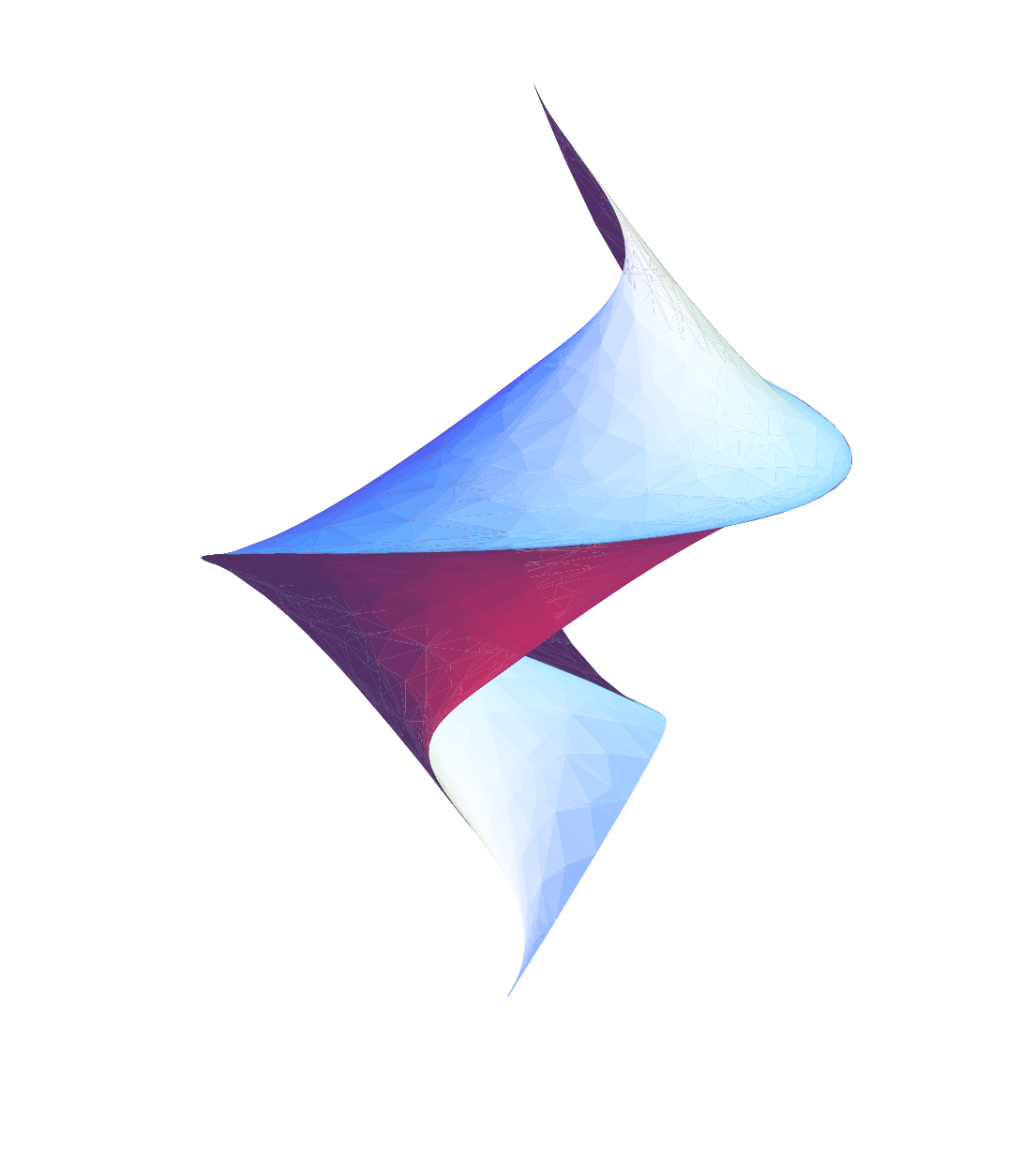}
\end{center}
\caption{Three swallowtails and cuspidal edges}\label{fi10}
\end{figure}

Following the same approach as in the non-characteristic case, see \cite{ACG, M2},  we can prove the existence of IAM with a prescribed curve of singularities (cuspidal edges or swallowtail) which determine the surface
 (see Figure \ref{fi10}).

\def\refname{References}

\end{document}